\documentclass[10 pt]{amsart}
\usepackage{amsfonts}
\usepackage{ifthen}
\usepackage{amsthm}
\usepackage{amsmath}
\usepackage{graphicx}
\usepackage{amscd,amssymb,amsthm}
\usepackage{graphicx}
\usepackage{epstopdf}
\usepackage{hyperref}

\newcounter{minutes}
\divide\time by 60
\newcounter{hours}
\multiply\time by 60 \addtocounter{minutes}{-\time}

\newtheorem{theorem}{Theorem}

\newcommand{\real}{\operatorname{Re}}

\keywords{Lommel, Struve and Bessel functions; univalent, starlike functions; radius of univalence and starlikeness;
zeros of Lommel, Struve and Bessel functions; Mittag-Leffler expansions; Laguerre-P\'olya class of entire functions.}
\subjclass[2010]{30C45, 30C15, 33C10}

\begin{document}

\title{Bounds for the radii of univalence of some special functions}

\author[I. Akta\c{s}]{Ibrahim Akta\c{s}}
\address{Department of Mathematical Engineering, Faculty of Engineering and Natural Sciences, G\"{u}m\"{u}\c{s}hane University, G\"{u}m\"{u}\c{s}hane 29100, Turkey}
\email{aktasibrahim38@gmail.com}
\author[\'A. Baricz]{\'Arp\'ad Baricz$^{\bigstar}$}
\address{Department of Economics, Babe\c{s}-Bolyai University, Cluj-Napoca
400591, Romania}
\address{Institute of Applied Mathematics, \'Obuda University, 1034
Budapest, Hungary}
\email{bariczocsi@yahoo.com}
\author[N. Ya\u{g}mur]{Nihat Ya\u{g}mur}
\address{Department of Mathematics, Erzincan University, Erzincan 24000,
Turkey}
\email{nhtyagmur@gmail.com}

\thanks{$^{\bigstar}$The research of \'A. Baricz was supported by the Romanian National Authority for Scientific Research, CNCS-UEFISCDI, under Grant
PN-II-RU-TE-2012-3-0190.}

\def\thefootnote{}
\footnotetext{ \texttt{File:~\jobname .tex,
          printed: \number\year-0\number\month-\number\day,
          \thehours.\ifnum\theminutes<10{0}\fi\theminutes}
} \makeatletter\def\thefootnote{\@arabic\c@footnote}\makeatother

\maketitle

\begin{abstract}
Tight lower and upper bounds for the radius of univalence of some normalized Bessel, Struve and Lommel functions of the
first kind are obtained via Euler-Rayleigh inequalities. It is shown also that the radius of univalence of the Struve
functions is greater than the corresponding radius of univalence of Bessel functions. Moreover, by using the idea of
Kreyszig and Todd, and Wilf it is proved that the radii of univalence of some normalized Struve and Lommel
functions are exactly the radii of starlikeness of the same functions. The Laguerre-P\'olya class of
entire functions plays an important role in our study.
\end{abstract}

\section{Introduction}

The surprising use of generalized hypergeometric functions in the solution of the famous Bieberbach conjecture created a considerable interest on geometric properties of special functions. Although in the sixties appeared only a small number of studies concerning geometric properties of Bessel functions (see \cite{brown,brown2,brown3,hayden,todd,wilf} and the references therein), in the last two decades the subject was reconsidered and many results were proved for Bessel, Struve and Lommel functions (see for example \cite{mathematica,publ,lecture,bdoy,bsk,bos,samy,basz,basz2,L1,szasz,szasz2} and the references therein).
The determination of the radii of starlikeness and convexity of Bessel, Struve and Lommel functions was intensively studied in the above papers, and it was shown that these radii are actually solutions of some transcendental equations. Since these solutions can be handled precisely only by computer softwares, it is natural to ask if there are some natural estimates for these radii. Motivated by this question, in this paper our aim is to present some tight lower and upper bounds for the radii of univalence (and starlikeness) of six normalized Bessel, Struve and Lommel functions of the first kind. Moreover, we show that these radii of univalence actually correspond to the radii of starlikeness. As in \cite{wilf}, the basic idea is that whereas the radius of univalence is quite troublesome to deal with directly, the radius of starlikeness is obtainable almost immediately from Weierstrassian's factorization. The results presented in this paper complement and improve the main results of the papers \cite{todd,wilf} and the proofs use the beautiful ideas from \cite{todd,wilf}. It is also important to mention that the method used in this paper is not new, its origins goes back to Euler and Rayleigh, see \cite{ismail} for more details. The bounds deduced for the
radii of univalence (and starlikeness) are actually particular cases of the well-known Euler-Rayleigh inequalities and it is possible to show that in the main results of this paper
the deduced lower bounds increase and the upper bounds decrease to the corresponding radii of univalence (and starlikeness). The fact that the lower bounds decrease can be
deduced directly from the corresponding Euler-Rayleigh inequalities, while the fact that the upper bounds are decreasing is actually a consequence of the Cauchy-Schwarz inequality. In other words, the inequalities
presented in this paper can be improved by using higher order Euler-Rayleigh inequalities. We restricted ourselves to the third Euler-Rayleigh inequalities since these
are already quite complicated. The organization of the paper is as follows: the next section contains the main results of the paper together with some remarks concerning the results presented therein, while the third
section is devoted to the proofs of the main results.

\section{Radii of univalence (and starlikeness) of some special functions}
\setcounter{equation}{0}

Let $\mathbb{D}_r=\{z\in\mathbb{C}:|z|<r\},$ where $r>0$ and let $\mathbb{D}_1=\mathbb{D}.$ Further let $f(z)=z+\alpha_2z^2+{\dots}$ be a univalent function.
We say that the radius of univalence of the function $f$ is the largest radius $r$ for which $f$ maps univalently the open disk $\mathbb{D}_r$ into some domain.
Similarly, the radius of starlikeness of the function $f$ is the largest radius $r$ for which $f$ maps $\mathbb{D}_r$ into a starlike domain with respect to origin.
In this paper we are interested on the radii of univalence (and starlikeness) of some normalized Bessel, Struve and Lommel functions of the first kind. We consider
the Bessel function of the first kind $J_{\nu}$, the Struve function of the first kind $\mathbf{H}_{\nu}$, and the Lommel
function of the first kind $s_{\mu,\nu}.$ The Bessel function of the first kind is a particular solution of the homogeneous Bessel differential equation, while the Struve
and Lommel functions of the first kind are particular solutions of the inhomogeneous Struve and Lommel differential equations. For more details on these functions we refer to
\cite{nist}. In the proof of the main results we are also using a special class of real entire functions, see \cite{DR} for more details. A real entire function $q$ belongs to the
Laguerre-P\'{o}lya class $\mathcal{LP}$ if it can be represented in the form
$$
q(x) = c x^{m} e^{-\alpha x^{2} + \beta x} \prod_{n\geq1}
\left(1+\frac{x}{x_{n}}\right) e^{-\frac{x}{x_{n}}},
$$
where $c,$ $\beta,$ $x_{n}$ are real numbers, $\alpha \geq 0,$ $m$ is a natural number or zero, and $\sum\limits_{n\geq 1} x_{n}^{-2}$ converges.

\subsection{Radii of univalence (and starlikeness) of Bessel functions of the first kind} The first main results we establish concerning the Bessel
function of the first kind read as follows.

\begin{theorem}\label{theoBes}
The radius of univalence (and starlikeness) $r^{\star}(\varphi_{\nu})$ of the normalized Bessel function of the first kind $z\mapsto \varphi_{\nu}(z)=2^{\nu}\Gamma(\nu+1)z^{1-\nu}J_{\nu}(z)$ satisfies $r^{\star}(\varphi_{\nu})<\sqrt{2(\nu+1)}$ for each $\nu>-1.$ Moreover, if $\nu>-1,$ then the radius of univalence satisfies
$$2\sqrt{\frac{\nu+1}{3}}<r^{\star}(\varphi_{\nu})<2\sqrt{\frac{3(\nu+1)(\nu+2)}{4\nu+13}},$$
$$2\sqrt[4]{\frac{(\nu+1)^2(\nu+2)}{4\nu+13}}<r^{\star}(\varphi_{\nu})<2\sqrt{\frac{(\nu+1)(\nu+3)(4\nu+13)}{2(4\nu^2+26\nu+49)}},$$
$$2\sqrt[6]{\frac{(\nu+1)^3(\nu+2)(\nu+3)}{2(4\nu^2+26\nu+49)}}<r^{\star}(\varphi_{\nu})<2\sqrt{\frac{2(\nu+1)(\nu+2)(\nu+4)(4\nu^2+26\nu+49)}{16\nu^4+208\nu^3+1032\nu^2+2341\nu+1987}}.$$
\end{theorem}

It is worth to mention that the inequality $r^{\star}(\varphi_{\nu})<\sqrt{2(\nu+1)},$ where $\nu>-1,$
improves the inequality $r^{\star}(\varphi_{\nu})<\sqrt{12(\nu+2)/5},$ obtained by Kreyszig and Todd \cite{todd}.
Note that our approach is more simple and is based on the first Rayleigh sum for the zeros of Bessel functions
of the first kind. We also note that the above first lower bound has been obtained also by
Kreyszig and Todd \cite{todd} by using a completely different method. Finally, we mention that if we consider the first
upper bound and the second lower bound for the radius of univalence $r^{\star}(\varphi_{\nu})$ from above, then we have
$$4\left(1+\frac{1}{\nu}\right)\sqrt{\frac{\nu+2}{4\nu+13}}<\frac{(r^{\star}(\varphi_{\nu}))^2}{\nu}<2\left(1+\frac{1}{\nu}\right),$$
which imply that $(r^{\star}(\varphi_{\nu}))^2/\nu\to 2$ as $\nu\to\infty.$ The above limit relation has been proved
earlier by Hayden and Merkes \cite{hayden} by using continued fractions for Bessel functions.

The next result is analogous to Theorem \ref{theoBes} and complements the main results in \cite{bsk,bos}.

\begin{theorem}\label{theoBes1}
Let $\nu>-1.$ The following assertions are true:
\begin{enumerate}
\item[\bf a.] The radius of univalence $r^{\star}(\phi_{\nu})$ of the normalized Bessel function of the first kind $z\mapsto \phi_{\nu}(z)=2^{\nu}\Gamma(\nu+1)z^{1-\frac{\nu}{2}}J_{\nu}(\sqrt{z})$ corresponds to the radius of starlikeness and it is the smallest positive root of the equation $\sqrt{z}J_{\nu}'(\sqrt{z})+(2-\nu)J_{\nu}(\sqrt{z})=0.$
\item[\bf b.] The radius of univalence (and starlikeness) $r^{\star}(\phi_{\nu})$ satisfies $r^{\star}(\phi_{\nu})<4(\nu+1)$ for each $\nu>-1.$ Moreover, if $\nu>-1,$ then the radius of univalence satisfies
$$2(\nu+1)<r^{\star}(\phi_{\nu})<\frac{8(\nu+1)(\nu+2)}{\nu+5},$$
$$\frac{4(\nu+1)\sqrt{\nu+2}}{\sqrt{\nu+5}}<r^{\star}(\phi_{\nu})<\frac{4(\nu+1)(\nu+3)(\nu+5)}{\nu^2+8\nu+23},$$
$$4(\nu+1)\sqrt[3]{\frac{(\nu+2)(\nu+3)}{\nu^2+8\nu+23}}<r^{\star}(\phi_{\nu})<\frac{4(\nu+1)(\nu+2)(\nu+4)(\nu^2+8\nu+23)}{\nu^4+15\nu^3+90\nu^2+267\nu+287}.$$
\item[\bf c.] The radius of univalence (and starlikeness) $r^{\star}(\phi_{\nu})$ satisfies $\frac{r^{\star}(\phi_{\nu})}{4(\nu+1)}\to 1$ as $\nu\to\infty.$ Moreover, we have the following asymptotic relation
$$r^{\star}(\phi_{\nu})=4(\nu+1)\left(1-\frac{1}{\nu}+\mathcal{O}\left(\nu^{-2}\right)\right),\ \ \ \mbox{as}\ \ \nu\to\infty.$$
\item[\bf d.] The function $\nu\mapsto r^{\star}(\phi_{\nu})$ is increasing on $(-1,\infty).$
\end{enumerate}
\end{theorem}

\subsection{Radii of univalence (and starlikeness) of Struve functions of the first kind} The next result concerning Struve functions of the first kind
complements the results in \cite{bdoy}.

\begin{theorem}\label{theo1} Let $\nu\in\left[-\frac{1}{2},\frac{1}{2}\right].$ The following assertions are true:
\begin{enumerate}
\item[\bf a.] The radius of univalence $r^{\star}(v_{\nu})$ of the normalized Struve function
$$z\mapsto v_{\nu}(z)=\sqrt{\pi }2^{\nu }z^{-\nu }\Gamma \left( \nu +\frac{3}{2}\right) \mathbf{H}_{\nu }(z)$$
corresponds to its radius of starlikeness, and it is the smallest positive root of the transcendental equation $z\mathbf{H}_{\nu}'(z)-\nu\mathbf{H}_{\nu}(z)=0.$
\item[\bf b.] The radius of univalence $r^{\star}(v_{\nu})$ satisfies $r^{\star}(v_{\nu})<\sqrt{3\left(\nu+\frac{3}{2}\right)}$ for $\nu\in\left[-\frac{1}{2},\frac{1}{2}\right].$ Moreover, when $\nu\in\left[-\frac{1}{2},\frac{1}{2}\right]$ the radius of univalence $r^{\star}(v_{\nu})$ satisfies the inequalities
    $$\sqrt{2\nu+3}<r^{\star}(v_{\nu})<\sqrt{\frac{3(2\nu+3)(2\nu+5)}{2\nu+9}},$$
    $$\sqrt[4]{\frac{3(2\nu+3)^2(2\nu+5)}{2\nu+9}}<r^{\star}(v_{\nu})<\sqrt{\frac{5(2\nu+3)(2\nu+7)(2\nu+9)}{3(4\nu^2+32\nu+79)}},$$
    $$\sqrt[6]{\frac{5(2\nu+3)^3(2\nu+5)(2\nu+7)}{4\nu^2+32\nu+79}}<r^{\star}(v_{\nu})<\sqrt{\frac{63(2\nu+3)(2\nu+5)(2\nu+9)(4\nu^2+32\nu+79)}{592\nu^4+9296\nu^3+56352\nu^2+159660\nu+171315}}.$$
\end{enumerate}
\end{theorem}

It is worth to mention that by using the inequalities $r^{\star}(v_{\nu})>\sqrt{2\nu+3}$ and $\sqrt{2(\nu+1)}>r^{\star}(\varphi_{\nu})$
it is clear that $r^{\star}(v_{\nu})>r^{\star}(\varphi_{\nu})$ when $\nu\in\left[-\frac{1}{2},\frac{1}{2}\right].$ Some numerical experiments
suggest that we can say more, that is, the positive zeros of the functions $z\mapsto z\mathbf{H}_{\nu}'(z)-\nu\mathbf{H}_{\nu}(z)$ and
$z\mapsto zJ_{\nu}'(z)-\nu J_{\nu}(z)$ are interlacing. However, we were unable to prove this result, which can be of interest for future research.

We also mention that by using the second lower bound for $r^{\star}(v_{\nu})$ in Theorem \ref{theo1} it can be shown that $r^{\star}(v_{\nu})>r^{\star}({v_{-\frac{1}{2}}})=\frac{\pi}{2}$ for each $\nu\in\left[\nu^{*},\frac{1}{2}\right],$ where $\nu^{*}\simeq-0.4935034122{\dots}$ is the unique root of the equation
$$\sqrt[4]{\frac{3(2\nu+3)^2(2\nu+5)}{2\nu+9}}=\frac{\pi}{2}.$$ Moreover, by using the third upper bound in part {\bf b} of Theorem \ref{theo1} we get that for $\nu\in\left[-\frac{1}{2},\frac{1}{2}\right]$ the radius of univalence $r^{\star}(v_{\nu})$ satisfies $r^{\star}(v_{\nu})\leq r^{\star}(v_{\frac{1}{2}}),$ where $r^{\star}(v_{\frac{1}{2}})\simeq2.33112237{\dots}$ is the smallest positive root of the trigonometric equation $z\sin z=1-\cos z.$ Taking into account these results it is natural to ask whether the function $\nu\mapsto r^{\star}(v_{\nu})$ is increasing on $\left[-\frac{1}{2},\frac{1}{2}\right].$ Note that if it would be possible to show that $\nu\mapsto h_{\nu,n}$ is increasing on $\left[-\frac{1}{2},\frac{1}{2}\right]$ for each fixed $n\in\mathbb{N},$ then the above result would be an immediate consequence of this. Here $h_{\nu,n}$ stands for the $n$th positive zero of the Struve function $\mathbf{H}_{\nu}.$

The following result complements the picture in \cite{bdoy} about Struve functions.

\begin{theorem}\label{theo2} Let $\nu\in\left[-\frac{1}{2},\frac{1}{2}\right].$ The following assertions are true:
\begin{enumerate}
\item[\bf a.] The radius of univalence $r^{\star}(w_{\nu})$ of the normalized Struve function $$z\mapsto w_{\nu}(z)=\sqrt{\pi }2^{\nu }z^{\frac{1-\nu}{2}}\Gamma \left( \nu +\frac{3}{2}\right) \mathbf{H}_{\nu }(\sqrt{z})$$ corresponds to its radius of starlikeness, and it is the smallest positive root of the transcendental equation $\sqrt{z}\mathbf{H}_{\nu}'(\sqrt{z})-(\nu-1)\mathbf{H}_{\nu}(\sqrt{z})=0.$
\item[\bf b.] The radius of univalence $r^{\star}(w_{\nu})$ satisfies $r^{\star}(w_{\nu})<3(2\nu+3)$ for $\nu\in\left[-\frac{1}{2},\frac{1}{2}\right].$ Moreover, when $\nu\in\left[-\frac{1}{2},\frac{1}{2}\right]$ the radius of univalence $r^{\star}(w_{\nu})$ satisfies the inequalities
$$\frac{3\left( 2\nu +3\right) }{2}<r^{\star}(w_{\nu})<\frac{15\left(2\nu +3\right) \left( 2\nu +5\right) }{2\nu +23},$$
$$\frac{3\left( 2\nu +3\right) \sqrt{5\left( 2\nu +5\right) }}{\sqrt{2\left(2\nu +23\right) }}<r^{\star}(w_{\nu})<\frac{21\left( 2\nu +3\right)
\left( 2\nu +7\right) \left( 2\nu +23\right) }{20\nu ^{2}+228\nu +1417},$$
$$3\left( 2\nu +3\right) \sqrt[3]{\frac{35\left( 2\nu +5\right) \left(2\nu+7\right) }{2\left( 20\nu ^{2}+228\nu +1417\right)}}<r^{\star}(w_{\nu})<\frac{15\left( 2\nu +3\right) \left( 2\nu +5\right) \left(2\nu +9\right) \left( 20\nu ^{2}+228\nu +1417\right) }{272\nu^{4}+5552\nu ^{3}+47584\nu ^{2}+247828\nu +416439}.$$
\end{enumerate}
\end{theorem}

We would like to mention that the starlikeness of the function $z\mapsto w_{\nu}(z)=\sqrt{z}v_{\nu}(\sqrt{z})$ can be obtained by using the interlacing property of the zeros Bessel and Struve functions of the first kind. More precisely, we have that if $|\nu|\leq\frac{1}{2}$ and $|z|=r<h_{\nu,1},$ then (see \cite{bdoy})
$$\real \frac{zw_{\nu}'(z)}{w_{\nu}(z)}=1-\real\sum_{n\geq 1}\frac{z}{h_{\nu,n}^2-z}\geq1-\sum_{n\geq 1}\frac{r}{h_{\nu,n}^2-r}=\frac{rw_{\nu}'(r)}{w_{\nu}(r)}>\frac{w_{\nu}'(1)}{w_{\nu}(1)}.$$
Here we used that $r\mapsto r/(h_{\nu,n}^2-r)$ is strictly increasing on $(0,h_{\nu,1})$ for each $n\in\mathbb{N}$ and consequently $r\mapsto {rw_{\nu}'(r)}/{w_{\nu}(r)}$ is strictly decreasing on $(0,h_{\nu,1}).$ Now, denoting by $j_{\nu,n}$ the $n$th positive zero of the Bessel function $J_{\nu}$ and by using the fact that the zeros of the Struve and Bessel functions of the first kind are interlacing, according to Steinig \cite{steinig}, it follows that for all $|\nu|\leq\frac{1}{2}$ we have
$$\frac{w_{\nu}'(1)}{w_{\nu}(1)}=1-\sum_{n\geq 1}\frac{1}{h_{\nu,n}^2-1}\geq 1-\sum_{n\geq 1}\frac{1}{j_{\nu,n}^2-1}>0,$$
where the last inequality follows from the proof of the starlikeness of the normalized Bessel functions of the first kind
$z\mapsto \phi_{\nu}(z)=2^{\nu}\Gamma(\nu+1)z^{1-\frac{\nu}{2}}J_{\nu}(\sqrt{z})$ and it is valid for
all $\nu\geq\nu_0,$ $\nu_0\simeq-0.5623$ being the root of the equation $\phi_{\nu}'(1)=0,$ see \cite{szasz} for more details. Thus, combining the above two chain of inequalities we obtain that indeed when $|\nu|\leq \frac{1}{2}$ we have that $w_{\nu}$ is starlike in the open unit disk. This result was proved recently by Baricz and Sz\'asz \cite{basz2}, but by using a somewhat different approach.

Finally, we mention that by using the inequalities $r^{\star}(w_{\nu})>\frac{3(2\nu+3)}{2}$ and $4(\nu+1)>r^{\star}(\phi_{\nu})$ we obtain that
$r^{\star}(w_{\nu})>r^{\star}(\phi_{\nu})$ for $\nu\in\left[-\frac{1}{2},\frac{1}{2}\right],$ that is, the radius of univalence of the normalized Struve function $w_{\nu}$ is greater than the corresponding radius of univalence of the normalized Bessel function $\phi_{\nu}.$

\subsection{Radii of univalence (and starlikeness) of Lommel functions of the first kind} The following results are analogous to the results obtained for Bessel
and Struve functions of the first kind and are related to the main results of \cite{bdoy} concerning Lommel functions of the first kind.

\begin{theorem}\label{theoLom1}
Let $\mu \in (-1,1),$ $\mu\neq0.$ The following assertions are true:
\begin{enumerate}
\item[\bf a.] The radius of univalence $r^{\star}(g_{\mu})$ of the normalized Lommel
function%
\[
z\mapsto g_{\mu }(z)=g_{\mu -\frac{1}{2},\frac{1}{2}}(z)=\mu (\mu +1)z^{-\mu
+\frac{1}{2}}s_{\mu -\frac{1}{2},\frac{1}{2}}(z)
\]
corresponds to its radius of starlikeness, and it is the smallest positive
root of the transcendental equation $zs_{\mu -\frac{1}{2},\frac{1}{2}%
}^{\prime }(z)-\left( \mu -\frac{1}{2}\right) s_{\mu -\frac{1}{2},\frac{1}{2}%
}(z)=0.$

\item[\bf b.] The radius of univalence $r^{\star}(g_{\mu})$ satisfies $r^{\star}(g_{\mu})<\sqrt{\frac{(\mu +2)(\mu +3)}{2}}.$
Moreover, the radius of univalence $r^{\star}(g_{\mu})$ satisfies the inequalities
\[
\sqrt{\frac{\left( \mu +2\right) \left( \mu +3\right) }{3}}<r^{\star}(g_{\mu})<%
\sqrt{\frac{3\left( \mu +2\right) \left( \mu +3\right) \left( \mu +4\right)
\left( \mu +5\right) }{-\mu ^{2}+31\mu +120}},
\]

\[
\sqrt[4]{\frac{\left( \mu +2\right) ^{2}\left( \mu +3\right) ^{2}\left( \mu
+4\right) \left( \mu +5\right) }{-\mu ^{2}+31\mu +120}}<r^{\star}(g_{\mu})<\sqrt{\frac{\left( \mu +2\right) \left( \mu +3\right) \left( \mu +6\right)
\left( \mu +7\right) \left( -\mu ^{2}+31\mu +120\right) }{3\left( \mu
^{4}-2\mu ^{3}+175\mu ^{2}+1842\mu +4032\right) }},
\]

\[
\sqrt[6]{\frac{\left( \mu +2\right) ^{3}\left( \mu +3\right) ^{3}\left( \mu
+4\right) \left( \mu +5\right) \left( \mu +6\right) \left( \mu +7\right) }{%
3\left( \mu ^{4}-2\mu ^{3}+175\mu ^{2}+1842\mu +4032\right) }}<r^{\star}(g_{\mu})
\]

\[
<\sqrt{\frac{3\left( \mu +2\right) \left( \mu +3\right) \left( \mu +4\right)
\left( \mu +5\right) \left( \mu +8\right) \left( \mu +9\right) \left( \mu
^{4}-2\mu ^{3}+175\mu ^{2}+1842\mu +4032\right) }{-\mu ^{8}+128\mu
^{7}+2196\mu ^{6}+24178\mu ^{5}+352645\mu ^{4}+3476958\mu
^{3}+17744328\mu ^{2}+44003088\mu +42301440}}.
\]
\end{enumerate}
\end{theorem}

The last main result is the following theorem concerning Lommel functions of the first kind.

\begin{theorem}\label{theoLom2}
Let $\mu \in (-1,1),$ $\mu\neq0.$ The following assertions are true:

\begin{enumerate}
\item[\bf a.] The radius of univalence $r^{\star}(h_{\mu})$ of the normalized Lommel
function
\[
z\mapsto h_{\mu }(z)=h_{\mu -\frac{1}{2},\frac{1}{2}}(z)=\mu (\mu +1)z^{%
\frac{3-2\mu }{4}}s_{\mu -\frac{1}{2},\frac{1}{2}}(\sqrt{z})
\]%
corresponds to its radius of starlikeness, and it is the smallest positive
root of the transcendental equation $2\sqrt{z}s_{\mu -\frac{1}{2},\frac{1}{2}%
}^{\prime }(\sqrt{z})-\left( 2\mu -3\right) s_{\mu -\frac{1}{2},\frac{1}{2}}(%
\sqrt{z})=0.$

\item[\bf b.] The radius of univalence $r^{\star}(h_{\mu})$ satisfies
$r^{\star}(h_{\mu})<(\mu +2)(\mu +3).$ Moreover, the radius of univalence $r^{\star}(h_{\mu})$ satisfies the inequalities
\[
\frac{\left( \mu +2\right) \left( \mu +3\right) }{2}<r^{\star}(h_{\mu})<\frac{%
\left( \mu +2\right) \left( \mu +3\right) \left( \mu +4\right) \left( \mu
+5\right) }{\left( -\mu ^{2}+3\mu +22\right) },
\]

\[
\frac{\left( \mu +2\right) \left( \mu +3\right) \sqrt{\left( \mu +4\right)
\left( \mu +5\right) }}{\sqrt{2\left( -\mu ^{2}+3\mu +22\right) }}<r^{\star}(h_{\mu})
<\frac{\left( \mu +2\right) \left( \mu +3\right) \left( \mu +6\right) \left(
\mu +7\right) \left( -\mu ^{2}+3\mu +22\right) }{\mu ^{4}-14\mu ^{3}-79\mu
^{2}+320\mu +1308},
\]

\[
\frac{\left( \mu +2\right) \left( \mu +3\right) \sqrt[3]{\left( \mu
+4\right) \left( \mu +5\right) \left( \mu +6\right) \left( \mu +7\right) }}{%
\sqrt[3]{2\left( \mu ^{4}-14\mu ^{3}-79\mu ^{2}+320\mu +1308\right) }}%
<r^{\star}(h_{\mu})
\]

\[
<\frac{-\left( \mu +2\right) \left( \mu +3\right) \left( \mu +4\right)
\left( \mu +5\right) \left( \mu +8\right) \left( \mu +9\right) \left( \mu
^{4}-14\mu ^{3}-79\mu ^{2}+320\mu +1308\right) }{3\mu ^{8}+96\mu
^{7}+2180\mu ^{6}+27338\mu ^{5}+164205\mu ^{4}+335662\mu
^{3}-881588\mu ^{2}-5070360\mu -6293376}.
\]
\end{enumerate}
\end{theorem}

\section{Proofs of the results}

\begin{proof}[\bf Proof of Theorem \ref{theoBes}] By using the first Rayleigh sum and the implicit relation for $r^{\star}(\varphi_{\nu}),$ obtained by Kreyszig and Todd \cite{todd}, we get for all $\nu>-1$ that
$$\frac{1}{(r^{\star}(\varphi_{\nu}))^{2}}=\sum_{n\geq 1}\frac{2}{j_{\nu,n}^2-(r^{\star}(\varphi_{\nu}))^2}>\sum_{n\geq 1}\frac{2}{j_{\nu,n}^2}=\frac{1}{2(\nu+1)}.$$

Now, by using the Euler-Rayleigh inequalities it is possible to have more tight bounds for the radius of univalence (and starlikeness) $r^{\star}(\varphi_{\nu}).$
For this recall that the zeros of $$\varphi_{\nu}(z)=\sum_{n\geq0}\frac{(-1)^nz^{2n+1}}{4^nn!(\nu+1)_n}$$ all are real when $\nu>-1.$
Consequently, this function belongs to the Laguerre-P\'olya class $\mathcal{LP}$ of real entire functions (see \cite{DC} for more details), which are uniform limits of
real polynomials whose all zeros are real. Now, since the Laguerre-P\'olya class $\mathcal{LP}$ is closed under differentiation, it follows
that $\varphi_{\nu}'$ belongs also to the Laguerre-P\'olya class and hence all of its zeros are real. Thus, the function
$z\mapsto \Psi_{\nu}(z)=\varphi_{\nu}'(2i\sqrt{z})$ has only negative real zeros and having growth order $\frac{1}{2}$ it can be written as the product
    $$\Psi_{\nu}(z)=\prod_{n\geq 1}\left(1+\frac{z}{a_{\nu,n}}\right),$$
    where $a_{\nu,n}>0$ for each $n\in\mathbb{N}.$ Now, by using the Euler-Rayleigh sum $\sigma_k=\sum_{n\geq1}a_{\nu,n}^{-k}$ and the infinite sum representation of the Bessel function $J_{\nu}$ we have (see the proof of \cite[Theorem 2]{wilf} for more details)
    $$\frac{\Psi_{\nu}'(z)}{\Psi_{\nu}(z)}=\sum_{n\geq1}\frac{1}{z+a_{\nu,n}}=\sum_{n\geq1}\sum_{k\geq0}\frac{(-1)^kz^k}{a_{\nu,n}^{k+1}}=\sum_{k\geq0}(-1)^k\sigma_{k+1}z^k,\ \ |z|<a_{\nu,1},$$
    $$\left.\frac{\Psi_{\nu}'(z)}{\Psi_{\nu}(z)}=\sum_{n\geq0}\frac{(2n+3)z^n}{n!(\nu+1)_{n+1}}\right/\sum_{n\geq0}\frac{(2n+1)z^n}{n!(\nu+1)_{n}}.$$
    From these relations it is possible to express the Euler-Rayleigh sums in terms of $\nu$ and by using the Euler-Rayleigh inequalities $\sigma_k^{-\frac{1}{k}}<a_{\nu,1}<\frac{\sigma_k}{\sigma_{k+1}}$ we obtain the inequalities for $\nu>-1$ and $k\in\mathbb{N}$
    $$2\sqrt{\sigma_k^{-\frac{1}{k}}}<r^{\star}(\varphi_{\nu})<2\sqrt{\frac{\sigma_k}{\sigma_{k+1}}}.$$
    Since $$\sigma_1=\frac{3}{\nu+1},\ \sigma_2=\frac{4\nu+13}{(\nu+1)^2(\nu+2)},\ \sigma_3=\frac{2(4\nu^2+26\nu+49)}{(\nu+1)^3(\nu+2)(\nu+3)}$$
    and $$\sigma_4=\frac{16\nu^4+208\nu^3+1032\nu^2+2341\nu+1987}{(\nu+1)^4(\nu+2)^2(\nu+3)(\nu+4)},$$
    in particular, when $k\in\{1,2,3\}$ from the above Euler-Rayleigh inequalities we have the next inequalities for $2\sqrt{a_{\nu,1}},$ that is,
    $$2\sqrt{\frac{\nu+1}{3}}<r^{\star}(\varphi_{\nu})<2\sqrt{\frac{3(\nu+1)(\nu+2)}{4\nu+13}},$$
    $$2\sqrt[4]{\frac{(\nu+1)^2(\nu+2)}{4\nu+13}}<r^{\star}(\varphi_{\nu})<2\sqrt{\frac{(\nu+1)(\nu+3)(4\nu+13)}{2(4\nu^2+26\nu+49)}},$$
    $$2\sqrt[6]{\frac{(\nu+1)^3(\nu+2)(\nu+3)}{2(4\nu^2+26\nu+49)}}<r^{\star}(\varphi_{\nu})<2\sqrt{\frac{2(\nu+1)(\nu+2)(\nu+4)(4\nu^2+26\nu+49)}{16\nu^4+208\nu^3+1032\nu^2+2341\nu+1987}},$$
and it is possible to have more tighter bounds for other values of $k\in\mathbb{N}.$
\end{proof}

\begin{proof}[\bf Proof of Theorem \ref{theoBes1}]
{\bf a.} We known that if the function $z\mapsto z+\alpha_2z^2+{\dots}$ has real coefficients, then its
radius of starlikeness is less or equal than its radius of univalence, see \cite{wilf}. On the other hand, we know that the radius of univalence
of the function $\phi_{\nu}$ is less or equal than the smallest positive zero of $\phi_{\nu}',$ according to Wilf \cite[p. 243]{wilf}. But,
the smallest positive zero of $\phi_{\nu}',$ that is, the first positive zero of the equation $\sqrt{z}J_{\nu}'(\sqrt{z})+(2-\nu)J_{\nu}(\sqrt{z})=0$ is actually the radius of starlikeness of $\phi_{\nu},$ according to \cite{bsk,bos}. These show that indeed
the radius of univalence corresponds to the radius of starlikeness of the function $\phi_{\nu}.$ Alternatively, we can follow Wilf's argument
(see the proof of \cite[Theorem 1]{wilf}) to show that the radii of univalence and starlikeness coincide.

{\bf b.} By using the well-known infinite product representation
$$2^{\nu}\Gamma(\nu+1)z^{-\nu}J_{\nu}(z)=\prod_{n\geq 1}\left(1-\frac{z^2}{j_{\nu,n}^2}\right)$$
we have that
\begin{equation}\label{logder}\frac{\phi_{\nu}'(z)}{\phi_{\nu}(z)}=\frac{1}{z}+\sum_{n\geq 1}\frac{1}{z-j_{\nu,n}^2},\end{equation}
which vanishes at $r^{\star}(\phi_{\nu}).$ In view of the first Rayleigh sum for the zeros of the Bessel functions of the first kind we get
\begin{equation}\label{radBes2}\frac{1}{r^{\star}(\phi_{\nu})}=\sum_{n\geq 1}\frac{1}{j_{\nu,n}^2-r^{\star}(\phi_{\nu})}>\sum_{n\geq 1}\frac{1}{j_{\nu,n}^2}=\frac{1}{4(\nu+1)}.\end{equation}

Now, consider the infinite sum representation of $\phi_{\nu}$ and its derivative%
\[
\phi_{\nu }(z)=\sum_{n\geq 0}\frac{\left( -1\right)
^{n}z^{n+1}}{n!4^{n}\left(\nu+1\right)_{n}},
\]%
\begin{equation}
\Phi_{\nu }(z)=\phi_{\nu }^{\prime }(-4z)=\sum_{n\geq
0}\frac{n+1}{\left(\nu+1\right)_{n}}\cdot\frac{z^n}{n!}.  \label{aBes}
\end{equation}%
The normalized Bessel function $\phi_{\nu }$ has only real zeros for $\nu >-1 $ and belongs to the Laguerre-P\'{o}lya class $\mathcal{LP}$ of real entire
functions, see \cite{DC}. Therefore $\phi_{\nu }^{\prime }$ belongs also to the Laguerre-P\'{o}lya class $\mathcal{LP}$ and has also only real zeros. Consequently, this
is also true for $\Phi_{\nu}.$ Moreover, since the coefficients of $\Phi_{\nu}(z)$ are non-negative and
$\Phi_{\nu}$ belongs to the Laguerre-P\'olya class $\mathcal{LP},$ it follows that $\Phi_{\nu}$ can have only negative zeros, and thus $\Phi_{\nu}(z)$ can be
written as the product
\begin{equation}\label{phirep}
\Phi_{\nu }(z)=\prod\limits_{n\geq 1}\left( 1+\frac{z}{b_{\nu ,n}}\right),\end{equation}
where $b_{\nu ,n}>0$ for each $n\in \mathbb{N}$. Now, by using the Euler-Rayleigh
sum $\rho_{k}=\sum_{n\geq 1}b_{\nu ,n}^{-k}$ and the infinite sum
representation of the Bessel function $J_{\nu }$ we have%
\[
\frac{\Phi_{\nu }^{\prime }(z)}{\Phi_{\nu }(z)}=\sum_{n\geq 1}\frac{1}{%
z+b_{\nu ,n}}=\sum_{n\geq 1}\sum_{k\geq 0}\frac{\left( -1\right) ^{k}z^{k}}{%
b_{\nu ,n}^{k+1}}=\sum_{k\geq 0}\left( -1\right) ^{k}\rho_{k+1}z^{k},\text{ \ \
}\left\vert z\right\vert <b_{\nu ,1},
\]%
\[
\frac{\Phi_{\nu }^{\prime }(z)}{\Phi_{\nu }(z)}=\left.\sum_{n\geq 0}\frac{( n+2)z^{n}}{n!\left( \nu+1\right)_{n+1}}\right/
\sum_{n\geq 0}\frac{(n+1)z^{n}}{n!\left(\nu+1\right)_{n}}.
\]
We can express the Euler-Rayleigh sums in terms of $\nu$ and by
using the Euler-Rayleigh inequalities $\rho_{k}^{-\frac{1}{k}}<b_{\nu ,1}<\frac{%
\rho_{k}}{\rho_{k+1}}$ we get the inequalities for $4b_{\nu ,1}$ when $\nu>-1$
and $k\in \mathbb{N}$%
\begin{equation}\label{euk23}4\rho_{k}^{-\frac{1}{k}}<r^{\star}(\phi_{\nu})<4\frac{\rho_{k}}{\rho_{k+1}}.\end{equation}
Since
\[
\rho_{1}=\frac{2}{\nu+1},\text{ \ }\rho_{2}=
\frac{\nu+5}{(\nu+1)^2(\nu+2)},\text{ \ }\rho_{3}=\frac{\nu^2+8\nu+23}{(\nu+1)^3(\nu+2)(\nu+3)}
\]%
and%
\[
\rho_{4}=\frac{\nu^4+15\nu^3+90\nu^2+267\nu+287}{(\nu+1)^4(\nu+2)^2(\nu+3)(\nu+4)},
\]%
in particular, when $k\in \{1,2,3\}$ we have the inequalities of this theorem.

{\bf c.} By using the Euler-Rayleigh inequality \eqref{euk23} for $k=2$ or $k=3$ we have that indeed the radius of univalence (and starlikeness) $r^{\star}(\phi_{\nu})$ satisfies $\frac{r^{\star}(\phi_{\nu})}{4(\nu+1)}\to 1$ as $\nu\to\infty.$ Moreover, if we use the inequality \eqref{euk23} for $k=3$ and the next asymptotic equalities
$4\rho_3^{-\frac{1}{3}}=4\frac{\rho_3}{\rho_4}=4(\nu+1)\left(1-\frac{1}{\nu}+\mathcal{O}\left(\nu^{-2}\right)\right)$ as $\nu\to\infty,$ we have the following asymptotic relation $r^{\star}(\phi_{\nu})=4(\nu+1)\left(1-\frac{1}{\nu}+\mathcal{O}\left(\nu^{-2}\right)\right)$ as $\nu\to\infty.$

{\bf d.} We denote the logarithmic derivative of $\phi_{\nu}$ by $\Delta_{\nu}.$ From \eqref{logder} and the left-hand side of \eqref{radBes2} it follows that $\Delta_{\nu}(z)>0$ for $z\in(0,r^{\star}(\phi_{\nu}))$ and $\Delta_{\nu}(r^{\star}(\phi_{\nu}))=0.$ Now, since the function $\nu\mapsto j_{\nu,n}$ is increasing on $(-1,\infty)$ for each fixed $n\in\mathbb{N},$ it follows that for $\mu>\nu$ and $z\in(0,r^{\star}(\phi_{\nu}))$ the terms of the series in
$$\Delta_{\mu}(z)-\Delta_{\nu}(z)=\sum_{n\geq 1}\frac{j_{\mu,n}^2-j_{\nu,n}^2}{(z-j_{\mu,n}^2)(z-j_{\nu,n}^2)}$$
are positive. Hence $\Delta_{\mu}(z)>0$ for $z\in(0,r^{\star}(\phi_{\nu}))$ and therefore $r^{\star}(\phi_{\mu})>r^{\star}(\phi_{\nu}).$
\end{proof}

\begin{proof}[\bf Proof of Theorem \ref{theo1}]
{\bf a.} The proof of this part concerning the radius of univalence goes along the lines introduced by Kreyszig and Todd \cite{todd}. First we show that the function $\theta\mapsto h(\theta)=\left|(re^{i\theta})^{-\nu}\mathbf{H}_{\nu}(re^{i\theta})\right|$ is increasing for $r<h_{\nu,1}$. By using the infinite product representation (see \cite{bps})
$$v_{\nu}(z)=z\prod_{n\geq 1}\left(1-\frac{z^2}{h_{\nu,n}^2}\right),$$
the expression $h(\theta)$ can be written as an infinite product of factors which are positive and increasing for $r<h_{\nu,1}.$ This shows that indeed $\theta\mapsto h(\theta)$ is increasing on $\left[0,\frac{\pi}{2}\right].$ On the other hand, if we consider $v_{\nu}(x)$ for $x\geq 0,$ then we have $v_{\nu}(x)=x+\mathcal{O}(x^3),$ and thus $v_{\nu}$ starts by increasing. Consequently, there is a last number $r^{\star}(v_{\nu})<h_{\nu,1}$ for which $v_{\nu}$ is a maximum. Moreover, the radius of univalence cannot exceed $r^{\star}(v_{\nu})$ because values of $v_{\nu}(x)$ are repeated for $x>r^{\star}(v_{\nu}).$ Now, by using again the above Weierstrassian decomposition of $v_{\nu}$ it is possible to show that $0\leq\arg (r_{\nu}e^{i\theta})^{-\nu}\mathbf{H}_{\nu}(r_{\nu}e^{i\theta})\leq\theta,$ where $r^{\star}(v_{\nu})$ is the smallest positive number for which $x\mapsto x^{-\nu}\mathbf{H}_{\nu}(x)$ has a maximum. So, if we consider the curve $\mathcal{C}$ consisting of three arcs: the segment $\mathcal{C}_1:$ $0\leq x\leq r^{\star}(v_{\nu})$ of the real axis; the arc $\mathcal{C}_2:$ $0<\theta<\frac{\pi}{2}$ of the circle $|z|=r^{\star}(v_{\nu})$; and the segment $\mathcal{C}_3:$ $r^{\star}(v_{\nu})>y>0$ of the imaginary axis; then by using the above properties it can be shown that if $\gamma_i$ is the map of $\mathcal{C}_i$ by $v_{\nu},$ then the arcs $\gamma_1,$ $\gamma_2$ and $\gamma_3$ are simple. Moreover, since $v_{\nu}$ is genuinely complex on $\mathcal{C}_2$ it follows that $\gamma_1$ or $\gamma_3$ cannot have common points with $\gamma_2.$ The arcs $\gamma_1$ and $\gamma_3$ cannot have also common points since $v_{\nu}$ is real on $\mathcal{C}_1$ and purely imaginary on $\mathcal{C}_3.$ Thus, the map $\gamma$ of $\mathcal{C}$ by $v_{\nu}$ has no double points and applying \cite[Lemma 4]{todd} it follows that the radius of univalence of $v_{\nu}$ is $r^{\star}(v_{\nu}).$

Now, according to Wilf \cite[Theorem 1]{wilf} we know that if the entire function $f$ has the form
$$f(z)=z\prod_{n\geq1}\left(1-\frac{z^2}{u_n^2}\right),$$
where $u_n>0$ and $\sum_{n\geq 1}u_n^{-2}<\infty,$ then the radii of univalence and starlikeness of $f$
coincide and are equal to the smallest positive zero of $f'.$  Thus, applying this result to the normalized
Struve function $v_{\nu}$ it follows that indeed $r^{\star}(v_{\nu})$ is the radius of univalence and starlikeness, which according to \cite{bdoy}
it is the smallest positive root of the transcendental equation $z\mathbf{H}_{\nu}'(z)-\nu\mathbf{H}_{\nu}(z)=0.$

{\bf b.} We proceed exactly as in the proof of Theorem \ref{theoBes} about the radius of univalence of normalized Bessel function discussed therein. By using the Mittag-Leffler expansion (see \cite{bps})
\begin{equation}\label{mle}\frac{\mathbf{H}_{\nu-1}(z)}{z\mathbf{H}_{\nu}(z)}-\frac{2\nu+1}{z^2}=\sum_{n\geq 1}\frac{2}{z^2-h_{\nu,n}^2},\end{equation}
the fact that
$$\frac{v_{\nu}'(z)}{v_{\nu}(z)}=-\frac{\nu}{z}+\frac{\mathbf{H}_{\nu}'(z)}{\mathbf{H}_{\nu}(z)}=-\frac{2\nu}{z}+\frac{\mathbf{H}_{\nu-1}(z)}{\mathbf{H}_{\nu}(z)}=
\frac{1}{z}+\sum_{n\geq 1}\frac{2z}{z^2-h_{\nu,n}^2}$$
vanishes at $r^{\star}(v_{\nu}),$ and taking the limit of both sides in \eqref{mle} as $z\to 0,$ we obtain the lower bound for the radius of univalence $r^{\star}(v_{\nu})$ in term of the first Rayleigh sum for the zeros of the Struve functions as follows
$$\frac{1}{(r^{\star}(v_{\nu}))^{2}}=\sum_{n\geq 1}\frac{2}{h_{\nu,n}^2-(r^{\star}(v_{\nu}))^2}>\sum_{n\geq 1}\frac{2}{h_{\nu,n}^2}=\frac{2}{3(2\nu+3)}.$$

Now, consider the infinite sum representation of $v_{\nu}$ and its derivative
$$v_{\nu}(z)=\frac{\sqrt{\pi}}{2}\sum_{n\geq0}\frac{(-1)^nz^{2n+1}}{4^n\Gamma\left(n+\frac{3}{2}\right)\left(\nu+\frac{3}{2}\right)_n},$$
$$V_{\nu}(z)=v_{\nu}'(2i\sqrt{z})=\frac{\sqrt{\pi}}{2}\sum_{n\geq0}\frac{(2n+1)z^{n}}{\Gamma\left(n+\frac{3}{2}\right)\left(\nu+\frac{3}{2}\right)_n}.$$
We know that $v_{\nu}$ has only real zeros for $\nu\in\left[-\frac{1}{2},\frac{1}{2}\right]$ and belongs to the Laguerre-P\'olya class $\mathcal{LP},$ see \cite{bps}.
Hence its derivative $v_{\nu}'$ has also only real zeros and thus $\varphi_{\nu}$ has only negative real zeros. Because the growth order
of the entire function $V_{\nu}$ is $\frac{1}{2}$ it follows that it can be written as the product
$$V_{\nu}(z)=\prod_{n\geq 1}\left(1+\frac{z}{c_{\nu,n}}\right),$$
where $c_{\nu,n}>0$ for each $n\in\mathbb{N}.$ Now, by using the Euler-Rayleigh sum $\tau_k=\sum_{n\geq1}c_{\nu,n}^{-k}$ and the
infinite sum representation of the Struve function $\mathbf{H}_{\nu}$ we have
$$\frac{V_{\nu}'(z)}{V_{\nu}(z)}=\sum_{n\geq1}\frac{1}{z+c_{\nu,n}}=\sum_{n\geq1}\sum_{k\geq0}\frac{(-1)^kz^k}{c_{\nu,n}^{k+1}}=\sum_{k\geq0}(-1)^k\tau_{k+1}z^k,\ \ |z|<c_{\nu,1},$$
$$\left.\frac{V_{\nu}'(z)}{V_{\nu}(z)}=\sum_{n\geq0}\frac{(2n+3)(n+1)z^n}{\Gamma\left(n+\frac{5}{2}\right)\left(\nu+\frac{3}{2}\right)_{n+1}}\right/
    \sum_{n\geq0}\frac{(2n+1)z^n}{\Gamma\left(n+\frac{3}{2}\right)\left(\nu+\frac{3}{2}\right)_{n}}.$$
Proceeding similarly as for Bessel functions in the proof of Theorem \ref{theoBes}, from the above relations we can express the Euler-Rayleigh sums in terms of $\nu$ and by using the Euler-Rayleigh inequalities $\tau_k^{-\frac{1}{k}}<c_{\nu,1}<\frac{\tau_k}{\tau_{k+1}}$ we obtain the inequalities for $2\sqrt{c_{\nu,1}}$ when $\nu\in\left[-\frac{1}{2},\frac{1}{2}\right]$ and $k\in\mathbb{N}$
$$2\sqrt{\tau_k^{-\frac{1}{k}}}<r^{\star}(v_{\nu})<2\sqrt{\frac{\tau_k}{\tau_{k+1}}}.$$
Since
$$\tau_1=\frac{4}{2\nu+3},\ \tau_2=\frac{16(2\nu+9)}{3(2\nu+3)^2(2\nu+5)},\ \tau_3=\frac{64(4\nu^2+32\nu+79)}{5(2\nu+3)^3(2\nu+5)(2\nu+7)}$$
and $$\tau_4=\frac{256(592\nu^4+9296\nu^3+56352\nu^2+159660\nu+171315)}{315(2\nu+3)^4(2\nu+5)^2(2\nu+7)(2\nu+9)},$$ in particular, when $k\in\{1,2,3\}$ we have the required inequalities.
\end{proof}

\begin{proof}[\bf Proof of Theorem \ref{theo2}]
{\bf a.} It is known that if the function $z\mapsto z+\alpha_2z^2+{\dots}$ has real coefficients, then its
radius of starlikeness is less or equal than its radius of univalence. On the other hand, we know that the radius of univalence
of the function $w_{\nu}$ is less or equal than the smallest positive zero of $w_{\nu}',$ according to Wilf \cite[p. 243]{wilf}. But,
the smallest positive zero of $w_{\nu}',$ that is, the first positive zero of the equation $\sqrt{z}\mathbf{H}_{\nu}'(\sqrt{z})+(1-\nu)\mathbf{H}_{\nu}(\sqrt{z})=0$ is actually the radius of starlikeness of $w_{\nu},$ according to \cite{bdoy}. These show that indeed
the radius of univalence corresponds to the radius of starlikeness of the function $w_{\nu}.$ Alternatively, we can follow Wilf's argument
(see the proof of \cite[Theorem 1]{wilf}) to show that the radii of univalence and starlikeness coincide.

{\bf b.} We proceed as in the proof of Theorem \ref{theoBes} about the radius of univalence of normalized Bessel function discussed therein. By using the infinite product representation (see \cite{bps})
$$\sqrt{\pi}2^{\nu}z^{-\nu-1}\Gamma\left(\nu+\frac{3}{2}\right)\mathbf{H}_{\nu}(z)=\prod_{n\geq 1}\left(1-\frac{z^2}{h_{\nu,n}^2}\right)$$
we have that
$$\frac{w_{\nu}'(z)}{w_{\nu}(z)}=\frac{1}{z}+\sum_{n\geq 1}\frac{1}{z-h_{\nu,n}^2},$$
which vanishes at $r^{\star}(w_{\nu}).$ In view of the first Rayleigh sum for the zeros of the Struve functions we get
$$\frac{1}{r^{\star}(w_{\nu})}=\sum_{n\geq 1}\frac{1}{h_{\nu,n}^2-r^{\star}(w_{\nu})}>\sum_{n\geq 1}\frac{1}{h_{\nu,n}^2}=\frac{1}{3(2\nu+3)}.$$

Now, consider the infinite sum representation of $w_{\nu}$ and its derivative%
\[
w_{\nu }(z)=\frac{\sqrt{\pi }}{2}\sum_{n\geq 0}\frac{\left( -1\right)
^{n}z^{n+1}}{4^{n}\Gamma \left( n+\frac{3}{2}\right) \left( \nu +\frac{3}{2}%
\right) _{n}},
\]%
\begin{equation}
W_{\nu }(z)=w_{\nu }^{\prime }(-4z)=\frac{\sqrt{\pi }}{2}\sum_{n\geq
0}\frac{\left( n+1\right)z^{n}}{\Gamma \left( n+\frac{3}{%
2}\right) \left( \nu +\frac{3}{2}\right) _{n}}=\frac{\sqrt{\pi }}{2}\sum_{n\geq
0}\frac{\left( n+1\right)!}{\Gamma \left( n+\frac{3}{2}\right) \left( \nu +\frac{3}{2}\right) _{n}}\cdot\frac{z^n}{n!}.  \label{a}
\end{equation}%
The normalized Struve function $w_{\nu }$ has only real zeros for $\nu \in \left[ -\frac{1}{2},\frac{1}{2}%
\right]$ and belongs to the Laguerre-P\'{o}lya class $\mathcal{LP}$ of real entire
functions, see \cite{bps}. Therefore $w_{\nu }^{\prime }$ belongs also to the Laguerre-P\'{o}lya class $\mathcal{LP}$ and has also only real zeros.
Consequently, this is also true for $W_{\nu}.$ Moreover, since the coefficients of $W_{\nu}(z)$ are non-negative and
$W_{\nu}$ belongs to the Laguerre-P\'olya class $\mathcal{LP},$ it follows that it has only negative real zeros, and thus, $W_{\nu}(z)$ can be
written as the product
\begin{equation}\label{psirep}
W_{\nu }(z)=\prod\limits_{n\geq 1}\left( 1+\frac{z}{d_{\nu ,n}}\right),\end{equation}
where $d_{\nu ,n}>0$ for each $n\in \mathbb{N}$. By using the Euler-Rayleigh
sum $\varsigma_{k}=\sum_{n\geq 1}d_{\nu ,n}^{-k}$ and the infinite sum
representation of the Struve function $\mathbf{H}_{\nu }$ we have%
\[
\frac{W_{\nu }^{\prime }(z)}{W_{\nu }(z)}=\sum_{n\geq 1}\frac{1}{%
z+d_{\nu ,n}}=\sum_{n\geq 1}\sum_{k\geq 0}\frac{\left( -1\right) ^{k}z^{k}}{%
d_{\nu ,n}^{k+1}}=\sum_{k\geq 0}\left( -1\right) ^{k}\varsigma_{k+1}z^{k},\text{ \ \
}\left\vert z\right\vert <d_{\nu ,1},
\]%
\[
\frac{W_{\nu }^{\prime }(z)}{W_{\nu }(z)}=\left.\sum_{n\geq 0}\frac{\left(
n+1\right) \left( n+2\right)z^{n}}{\Gamma \left( n+\frac{5}{2}%
\right) \left( \nu +\frac{3}{2}\right) _{n+1}}\right/ \sum_{n\geq 0}\frac{%
\left( n+1\right)z^{n}}{\Gamma \left( n+\frac{3}{2}\right) \left(
\nu +\frac{3}{2}\right) _{n}}.
\]%
It is possible to express the Euler-Rayleigh sums in terms of $\nu$ and by
using the Euler-Rayleigh inequalities $\varsigma_{k}^{-\frac{1}{k}}<d_{\nu ,1}<\frac{%
\varsigma_{k}}{\varsigma_{k+1}}$ we get the inequalities for $4d_{\nu ,1}$ when $\nu\in\left[-\frac{1}{2},\frac{1}{2}\right]$
and $k\in \mathbb{N}$%
$$4\varsigma_{k}^{-\frac{1}{k}}<r^{\star}(w_{\nu})<4\frac{\varsigma_{k}}{\varsigma_{k+1}}.$$
Since
\[
\varsigma_{1}=\frac{8}{3\left( 2\nu +3\right) },\text{ \ }\varsigma_{2}=
\frac{32\left( 2\nu +23\right) }{45\left( 2\nu +3\right) ^{2}\left(
2\nu +5\right) },\text{ \ }\varsigma_{3}=\frac{128\left( 20\nu ^{2}+228\nu
+1417\right) }{945\left( 2\nu +3\right) ^{3}\left( 2\nu +5\right) \left(
2\nu +7\right) }
\]%
and%
\[
\varsigma_{4}=\frac{512\left( 272\nu ^{4}+5552\nu ^{3}+47\,584\nu
^{2}+247\,828\nu +416\,439\right) }{14175\left( 2\nu +3\right) ^{4}\left( 2\nu
+5\right) ^{2}\left( 2\nu +7\right) \left( 2\nu +9\right) },
\]%
in particular, when $k\in \{1,2,3\}$ we have the required inequalities.
\end{proof}

\begin{proof}[\bf Proof of Theorem \ref{theoLom1}]
First we prove the results of this theorem for $\mu\in(0,1).$

{\bf a.} Recall that if the function $z\mapsto z+\alpha_2z^2+{\dots}$ has real coefficients, then its
radius of starlikeness is less or equal than its radius of univalence. On the other hand, we know that the radius of univalence
of the function $g_{\mu}$ is less or equal than the smallest positive zero of $g_{\mu}',$ according to Wilf \cite[p. 243]{wilf}. But,
the smallest positive zero of $g_{\mu}',$ that is, the first positive zero of the equation $zs_{\mu -\frac{1}{2},\frac{%
1}{2}}^{\prime }(z)-\left( \mu -\frac{1}{2}\right) s_{\mu -\frac{1}{2},\frac{1}{2}}(z)=0$ is actually the radius of starlikeness of $g_{\mu},$ according to \cite{bdoy}. These show that indeed the radius of univalence corresponds to the radius of starlikeness of the function $g_{\mu}.$ Alternatively, we can follow Wilf's argument (see the proof of \cite[Theorem 1]{wilf}) to show that the radii of univalence and starlikeness coincide.

{\bf b.} By using the infinite product representation (see \cite{koumandos,L1})
\begin{equation}\label{infLom}
s_{\mu -\frac{1}{2},\frac{1}{2}}(z)=\frac{z^{\mu +\frac{1}{2}}}{\mu (\mu +1)}%
\prod\limits_{n\geq 1}\left( 1-\frac{z^{2}}{\xi _{\mu ,n}^{2}}\right) ,
\end{equation}
we get%
\[
\frac{g_{\mu }^{\prime }(z)}{g_{\mu }(z)}=\frac{1}{z}+\sum_{n\geq 0}\frac{2z%
}{z^{2}-\xi _{\mu ,n}^{2}},
\]%
which vanishes at $r^{\star}(g_{\mu}).$ Here $\xi _{\mu ,n}$ stands for the $n$th positive
zero of the Lommel function $s_{\mu -\frac{1}{2},\frac{1}{2}}.$ In view of the first Rayleigh sum for the
zeros of the Lommel functions we obtain%
\[
\frac{1}{\left(r^{\star}(g_{\mu})\right) ^{2}}=\sum_{n\geq 0}\frac{2}{\xi _{\mu
,n}^{2}-\left(r^{\star}(g_{\mu})\right) ^{2}}>\sum_{n\geq 0}\frac{2}{\xi _{\mu
,n}^{2}}=\frac{2}{\left( \mu +2\right) \left( \mu +3\right) }.
\]

Now, consider the infinite sum representation of $g_{\mu}$ and its
derivative%
\[
g_{\mu }(z)=\sum_{n\geq 0}\frac{\left( -1\right) ^{n}z^{2n+1}}{4^{n}\left(
\frac{\mu +2}{2}\right) _{n}\left( \frac{\mu +3}{2}\right) _{n}},
\]%
\[
G_{\mu }(z)=g_{\mu }^{\prime }(2i\sqrt{z})=\sum_{n\geq 0}\frac{\left(
2n+1\right) z^{n}}{\left( \frac{\mu +2}{2}\right) _{n}\left( \frac{\mu +3}{2}%
\right) _{n}}.
\]%
The function $g_{\mu }$ has only real zeros for $\mu \in (0,1)$ and belongs to the
Laguerre-P\'{o}lya class $\mathcal{LP}$ of real entire functions, see \cite{koumandos}. Hence $g_{\mu }^{\prime }$ belongs also to the
Laguerre-P\'{o}lya class $\mathcal{LP}$. Thus, $G_{\mu }$ has only negative real zeros and it can be written as the product%
\[
G_{\mu }(z)=\prod\limits_{n\geq 1}\left( 1+\frac{z}{\kappa_{\mu ,n}}\right) ,
\]%
where $\kappa_{\mu ,n}>0$ for each $n\in \mathbb{N}$. By using the Euler-Rayleigh
sum $\varrho_{k}=\sum_{n\geq 1}\kappa_{\mu ,n}^{-k}$ and the infinite sum
representation of the Lommel function $s_{\mu -\frac{1}{2},\frac{1}{2}}$ we
have%
\[
\frac{G_{\mu }^{\prime }(z)}{G_{\mu }(z)}=\sum_{n\geq 1}\frac{1}{z+\kappa_{\mu ,n}%
}=\sum_{n\geq 1}\sum_{k\geq 0}\frac{\left( -1\right) ^{k}z^{k}}{\kappa_{\mu
,n}^{k+1}}=\sum_{k\geq 0}\left( -1\right) ^{k}\varrho_{k+1}z^{k},\text{ \ \ }%
\left\vert z\right\vert <\kappa_{\mu ,1},
\]%
\[
\frac{G_{\mu }^{\prime }(z)}{G_{\mu }(z)}=\left.\sum_{n\geq 0}\frac{\left(
n+1\right) \left( 2n+3\right) z^{n}}{\left( \frac{\mu +2}{2}\right)
_{n+1}\left( \frac{\mu +3}{2}\right) _{n+1}}\right/ \sum_{n\geq 0}\frac{%
\left( 2n+1\right) z^{n}}{\left( \frac{\mu +2}{2}\right) _{n}\left( \frac{%
\mu +3}{2}\right) _{n}}.
\]%
It is possible to express the Euler-Rayleigh sums in terms of $\mu $ and by
using the Euler-Rayleigh inequalities $\varrho_{k}^{-\frac{1}{k}}<\kappa_{\mu ,1}<\frac{%
\varrho_{k}}{\varrho_{k+1}}$ we get the inequalities for $2\sqrt{\kappa_{\mu ,1}}$ when $\mu
\in (0,1)$ and $k\in \mathbb{N}$%
\[
2\sqrt{\varrho_{k}^{-\frac{1}{k}}}<r^{\star}(g_{\mu})<2\sqrt{\frac{\varrho_{k}}{\varrho_{k+1}}}.
\]%
Since
\[
\varrho_{1}=\frac{12}{\left( \mu +2\right) \left( \mu +3\right) },\text{ \ }%
\varrho_{2}=\allowbreak \frac{16\left( -\mu ^{2}+31\mu +120\right) }{\left( \mu
+2\right) ^{2}\left( \mu +3\right) ^{2}\left( \mu +4\right) \left( \mu
+5\right) },
\]%
\[
\varrho_{3}=\allowbreak \frac{192\left( \mu ^{4}-2\mu ^{3}+175\mu ^{2}+1842\mu
+4032\right) }{\left( \mu +2\right) ^{3}\left( \mu +3\right) ^{3}\left( \mu
+4\right) \left( \mu +5\right) \left( \mu +6\right) \left( \mu +7\right) }
\]%
and%
\[
\varrho_{4}=\frac{-\mu ^{8}+128\mu ^{7}+2196\mu ^{6}+24178\mu ^{5}+352645\mu ^{4}+
3476958\mu ^{3}+17744328\mu ^{2}+44003088\mu +42301440 }{256^{-1}\left( \mu +2\right) ^{4}\left( \mu +3\right) ^{4}\left( \mu
+4\right) ^{2}\left( \mu +5\right) ^{2}\left( \mu +6\right) \left( \mu
+7\right) \left( \mu +8\right) \left( \mu +9\right) }.
\]%
in particular, when we take $k\in \{1,2,3\}$ in the above Euler-Rayleigh inequalities, we have the inequalities of this theorem.

Now, we prove that parts {\bf a} and {\bf b} also hold when $\mu
\in \left( -1,0\right).$ In order to do this, suppose that $\mu \in (0,1)$
and repeat the above proof, substituting $\mu $ by $\mu -1$, $g_{\mu}$ by
the function $g_{\mu-1}$ and taking into account that the $n$th positive
zero of $g_{\mu-1},$ denoted by $\zeta _{\mu,n},$ are all real, since $g_{\mu-1}$
belongs also to the Laguerre-P\'{o}lya class $\mathcal{LP}$ of entire functions (see \cite{koumandos}).
It is worth mentioning that
\[
\frac{g_{\mu -1}^{\prime }(z)}{g_{\mu -1}(z)}=\frac{1}{z}+\sum_{n\geq 0}%
\frac{2z}{z^{2}-\zeta _{\mu ,n}^{2}}
\]%
holds for $\mu\in\left(0,1\right)$. In this case we have%
\[
\sum_{n\geq 0}\frac{2}{\zeta _{\mu ,n}^{2}}=\frac{2}{\left( \mu +1\right)
\left( \mu +2\right) },
\]%
hence the parts {\bf a} and {\bf b} are valid for $\mu -1$ instead
of $\mu .$ Thus, now replacing $\mu $ by $\mu +1$, we obtain the statements of
the parts {\bf a} and {\bf b} for $\mu\in\left(-1,0\right).$
\end{proof}

\begin{proof}[\bf Proof of Theorem \ref{theoLom2}]
We proceed similarly as in the proof of Theorem \ref{theoLom1}. First we show that the statements of this theorem are valid for $\mu\in(0,1),$
and then by using a similar argument as in the proof of Theorem \ref{theoLom1} we consider the case of $\mu\in(-1,0).$ Since the argument concerning
the case when $\mu\in(-1,0)$ goes along the lines introduced at the end of the proof of Theorem \ref{theoLom1}, we omit the details.

{\bf a.} We shall use again the fact that if the function $z\mapsto z+\alpha_2z^2+{\dots}$ has real coefficients, then its
radius of starlikeness is less or equal than its radius of univalence. On the other hand, we know that the radius of univalence
of the function $h_{\mu}$ is less or equal than the smallest positive zero of $h_{\mu}',$ according to Wilf \cite[p. 243]{wilf}. But,
the smallest positive zero of $h_{\mu}',$ that is, the first positive zero of the equation $2\sqrt{z}s_{\mu -\frac{1}{2%
},\frac{1}{2}}^{\prime }(\sqrt{z})-\left( 2\mu -3\right) s_{\mu -\frac{1}{2},%
\frac{1}{2}}(\sqrt{z})=0$ is actually the radius of starlikeness of $h_{\mu},$ according to \cite{bdoy}. These show that indeed the radius of univalence corresponds to the radius of starlikeness of the function $h_{\mu}.$ Alternatively, we can follow Wilf's argument (see the proof of \cite[Theorem 1]{wilf}) to show that the radii of univalence and starlikeness coincide.

{\bf b.} By using the infinite product representation \eqref{infLom} we get
\[
\frac{h_{\mu }^{\prime }(z)}{h_{\mu }(z)}=\frac{1}{z}+\sum_{n\geq 0}\frac{1}{%
z-\xi _{\mu ,n}^{2}},
\]%
which vanishes at $r^{\star}(h_{\mu}).$ In view of the first Rayleigh sum for the
zeros of the Lommel functions we obtain%
\[
\frac{1}{r^{\star}(h_{\mu})}=\sum_{n\geq 0}\frac{1}{\xi _{\mu ,n}^{2}-r^{\star}(h_{\mu})}%
>\sum_{n\geq 0}\frac{1}{\xi _{\mu ,n}^{2}}=\frac{1}{\left( \mu +2\right)
\left( \mu +3\right) }.
\]

Now, consider the infinite sum representation of $h_{\mu }$ and its
derivative%
\[
h_{\mu }(z)=\sum_{n\geq 0}\frac{\left( -1\right) ^{n}z^{n+1}}{4^{n}\left(
\frac{\mu +2}{2}\right) _{n}\left( \frac{\mu +3}{2}\right) _{n}},
\]%
\[
H_{\mu }(z)=h_{\mu }^{\prime }(-4z)=\sum_{n\geq 0}\frac{\left( n+1\right)
z^{n}}{\left( \frac{\mu +2}{2}\right) _{n}\left( \frac{\mu +3}{2}\right) _{n}%
}.
\]%
Note that $h_{\mu}$ has only real zeros for $\mu \in (0,1)$ and belongs to the
Laguerre-P\'{o}lya class $\mathcal{LP}$ of real entire functions, see \cite{koumandos}. Hence $h_{\mu}^{\prime}$ belongs also to $\mathcal{LP}$ and
has also only real zeros. Since the coefficients of $H_{\mu }(z)$ are
non-negative and $H_{\mu }$\ belongs to the Laguerre-P\'olya class $\mathcal{LP%
}$, it has only negative zeros, and thus $H_{\mu }$ can
be written as the product%
\[
H_{\mu }(z)=\prod\limits_{n\geq 1}\left( 1+\frac{z}{\lambda _{\mu ,n}}%
\right) ,
\]%
where $\lambda _{\mu ,n}>0$ for each $n\in \mathbb{N}$. By using the
Euler-Rayleigh sum $\eta _{k}=\sum_{n\geq 0}\lambda _{\mu ,n}^{-k}$ and the
infinite sum representation of the Lommel function $s_{\mu -\frac{1}{2},%
\frac{1}{2}}$ we have%
\[
\frac{H_{\mu }^{\prime }(z)}{H_{\mu }(z)}=\sum_{n\geq 1}\frac{1}{z+\lambda
_{\mu ,n}}=\sum_{n\geq 1}\sum_{k\geq 0}\frac{\left( -1\right) ^{k}z^{k}}{%
\lambda _{\mu ,n}^{k+1}}=\sum_{k\geq 0}\left( -1\right) ^{k}\eta _{k+1}z^{k},%
\text{ \ \ }\left\vert z\right\vert <\lambda _{\mu ,1},
\]%
\[
\frac{H_{\mu }^{\prime }(z)}{H_{\mu }(z)}=\left.\sum_{n\geq 0}\frac{\left(
n+1\right) \left( n+2\right) z^{n}}{\left( \frac{\mu +2}{2}\right)
_{n+1}\left( \frac{\mu +3}{2}\right) _{n+1}}\right/ \sum_{n\geq 0}\frac{%
\left( n+1\right) z^{n}}{\left( \frac{\mu +2}{2}\right) _{n}\left( \frac{\mu
+3}{2}\right) _{n}}.
\]%
Similarly as in the above proofs, it is possible to express the Euler-Rayleigh sums in terms of $\mu $ and by
using the Euler-Rayleigh inequalities $\eta _{k}^{-\frac{1}{k}}<\lambda
_{\mu ,1}<\frac{\eta _{k}}{\eta _{k+1}}$ we get the inequalities for $%
4\lambda _{\mu ,1}$ when $\mu \in (0,1)$ and $k\in \mathbb{N}$%
\[
4\eta _{k}^{\frac{-1}{k}}<r_{\mu}(h)<\frac{4\eta _{k}}{\eta _{k+1}}.
\]%
Since
\[
\eta _{1}=\frac{8}{\left( \mu +2\right) \left( \mu +3\right) },\text{ \ }%
\eta _{2}=\allowbreak \frac{32\left( -\mu ^{2}+3\mu +22\right) }{\left( \mu
+2\right) ^{2}\left( \mu +3\right) ^{2}\left( \mu +4\right) \left( \mu
+5\right) },
\]%
\[
\eta _{3}=\allowbreak \frac{128\left( \mu ^{4}-14\mu ^{3}-79\mu ^{2}+320\mu
+1308\right) }{\left( \mu +2\right) ^{3}\left( \mu +3\right) ^{3}\left( \mu
+4\right) \left( \mu +5\right) \left( \mu +6\right) \left( \mu +7\right) }
\]%
and%
\[
\eta _{4}=\frac{-3\mu ^{8}-96\mu ^{7}-2180\mu ^{6}-27338\mu ^{5}-164205\mu ^{4}-335662\mu
^{3}+ 881588\mu ^{2}+5070360\mu +6293376}{512^{-1}\left( \mu +2\right) ^{4}\left( \mu +3\right) ^{4}\left( \mu
+4\right) ^{2}\left( \mu +5\right) ^{2}\left( \mu +6\right) \left( \mu
+7\right) \left( \mu +8\right) \left( \mu +9\right) },
\]%
in particular, when $k\in \{1,2,3\}$ we get the inequalities of the theorem.
\end{proof}

\end{document}